\documentclass{elsart}

\usepackage{amsfonts}
\usepackage{amstext}
\usepackage{amssymb,amsmath}

\def\N{\mathbb N}\def\R{\mathbb R}
\def\EE{\mathbb E}\def\PP{\mathbb P}
\def\0{\mathbf 0}
\def\1{\mathbf 1}

\def\q{\mathbf q}
\def\r{\mathbf r}
\def\x{\mathbf x}
\def\uu{\mathbf u}
\def\vv{\mathbf v}
\def\w{\mathbf w}
\def\z{\mathbf z}
\def\y{\mathbf y}
\def\X{\mathbf X}
\def\V{\mathbf V}
\def\Z{\mathbf Z}
\def\A{\mathbf A}
\def\Y{\mathbf Y}
\def\EEE{\mathbf E}
\def\PPP{\mathbf P}
\def\RRR{\mathbf R}
\def\W{\mathbf W}
\def\BB{\mathbf B}
\def\DD{\mathbf D}
\def\QQ{\mathbf Q}
\def\SS{\mathbf S}

\def\eps{\varepsilon}

\def\diag{\text{diag}\,}
\def\rk{\text{rank}\,}
\def\dist{\text{dist}}
\def\Span{\text{Span}\,}
\def\Var{{\text{Var}}\,}

\begin{document}

\begin{frontmatter}
\title { Singular value decomposition of large random matrices 
\\ (for two-way classification of microarrays)}

\author[bmemat]{Marianna Bolla\corauthref{ca}}, 
\ead{marib@math.bme.hu}
\corauth[ca]{Corresponding author.}
\author[bmeszit]{Katalin Friedl},
\ead{friedl@cs.bme.hu}
\author[szeged]{Andr\'as Kr\'amli}
\ead{kramli@informatika.ilab.sztaki.hu}

\address[bmemat]{Institute of Mathematics,
Budapest University of Technology and Economics}
\address[bmeszit]{Dept. of Computer Science,
Budapest University of Technology and Economics }
\address[szeged]{Bolyai Institute, University of Szeged}


\begin{abstract} 
Asymptotic behavior of the singular value decomposition (SVD)
of  blown up matrices and normalized blown up 
contingency tables  exposed  to Wigner-noise is  investigated.
It is proved that such an $m\times n$ matrix almost surely has a 
constant number of large singular values (of order $\sqrt{mn}$),
 while the rest of the singular values are of order $\sqrt{m+n}$
as $m,n\to\infty$.
Concentration results of Alon at al. for the eigenvalues of large symmetric
random matrices are adapted
to the rectangular case, and on this basis,
almost sure results for the singular values as well as for the
corresponding isotropic subspaces are proved. An algorithm, applicable to 
two-way classification of microarrays,
is also given that finds the underlying block structure.
\end{abstract}


\end{frontmatter}

\section{Introduction}

%

A general problem of multivariate statistics is to find linear structures
in large real-world data sets like  internet or microarray measurements. 
In~\cite{bolla2}, large symmetric blown up matrices burdened 
with a so-called symmetric Wigner-noise were investigated. 
It was proved that such an $n\times n$
matrix has some protruding eigenvalues (of order $n$), while the majority 
of the eigenvalues is at most of order $\sqrt{n}$ with probability
tending to 1  as $n\to\infty$. 
Our goal is to generalize these results for  the stability of
SVD of large rectangular random matrices
and to apply them for  the contingency table matrix formed by 
categorical variables in order to perform two-way clustering of these
variables.
First we introduce some notation.

\begin{defn}\label{wigner}
The  $m\times n$ real matrix $\W$ is a
{\em  Wigner-noise} if its
 entries $w_{ij}$ $(1\le i\le m, \, 1\le j\le n)$
are independent random variables, $\EE (w_{ij}) =0$,
and the $w_{ij}$'s  are uniformly bounded
(i.e., there is a constant $K>0$, independently of $m$ and $n$, 
 such that $|w_{ij}| \le K$, $\forall i,j$).
\end{defn}

Though, the main results of this paper can be extended to 
$w_{ij}$'s with any light-tail distribution (especially  to
Gaussian distributed $w_{ij}$'s), our almost sure results will be based
on the assumptions of Definition~\ref{wigner}.
%
%

\begin{defn}\label{blownup}
The $m\times n$ real matrix $\BB$  is a {\em blown up  matrix}, if there is an 
$a\times b$ so-called {\it pattern matrix} $\PPP$ with entries 
$0\le p_{ij} \le 1$, 
and there  are positive integers 
$m_1 ,\dots ,m_a$ with $\sum_{i=1}^a m_i =m$ and
$n_1 ,\dots ,n_b$ with $\sum_{i=1}^b n_i =n$,
such that
the matrix $\BB$ can be divided into $a\times b$ blocks, where 
block $(i,j)$ is 
an $m_i \times n_j$ matrix with entries 
equal to $p_{ij}$ $(1\le i\le a, \, 1\le j\le b )$.
\end{defn}

Such schemes are sought for in microarray analysis and 
they are called chess-board patterns, cf.~\cite{KBCG}.
Let us fix the matrix $\PPP$,  blow it up to obtain matrix $\BB$,
and let $\A =\BB +\W$, where $\W$ is a Wigner-noise of appropriate size.
We are interested in  the  properties of $\A$ when 
 $m_1 ,\dots ,m_a\to\infty$ and $n_1 ,\dots ,n_b
\to \infty$, roughly speaking, 
at the same rate.
More precisely,
we make two different constraints on the growth of the sizes
$m$, $n$, and the growth rate of their components. 
The first one is needed for  all our reasonings, while the second one
will be used in the case of noisy correspondence matrices, only.

\begin{defn}
~
\begin{description}
\item[GC1] (Growth Condition 1) \\
 There exists a constant $0<c<1$
such that $m_i  /m \ge c$  ($i=1,\dots ,a$) and 
 there exists a constant $0<d<1$ such that
     $n_i  /n \ge d$   ($i=1,\dots ,b$).
\item[GC2]  (Growth Condition 2) \\
 There exist constants $C\ge 1$, $D\ge 1$, and $C_0 >0$, $D_0 > 0$
such that 
\, $m\le C_0 \cdot n^C$ and $n\le D_0 \cdot m^D$ hold
 for sufficiently large $m$ and $n$.
\end{description}
\end{defn}

\begin{rem}
~
{\it GC1}  implies that 
\begin{equation}\label{const}
  c \le \frac{m_k}{m_i} \le \frac1{c} \quad\text{and}\quad
  d \le \frac{n_{\ell}}{n_j} \le \frac1{d}   
\end{equation}
hold for any pair of indices $k,i \in \{ 1,\dots ,a \}$ and
$\ell,j \in \{ 1,\dots ,b \}$.
\end{rem}

We want to establish some property ${\cal P}_{m,n}$ that holds for the 
$m\times n$ random matrix $\A =\BB +\W$ (briefly, $\A_{m\times n}$)
with $m$ and $n$ large enough.
In this paper  ${\cal P}_{m,n}$ is  mostly related to 
the SVD of $\A_{m\times n}$. 

\begin{defn}\label{conv}
Property ${\cal P}_{m,n}$ holds for $\A_{m\times n}$
almost surely (with probability 1) if  
$
\PP\left( \exists \,\, m_0, n_0 \in \N \,\, \text{such that for} \,\,
 m\ge m_0 \,\, \text{and} \,\, n\ge n_0 \,\, 
 \A_{m\times n} \,\, \text{has} \,\, {\cal P}_{m,n} \right)=1.$
Here we may assume {\it GC1} or {\it GC2} for the growth of $m$ and $n$,
while $K$ is kept fixed.
\end{defn}

In combinatorics literature convergence in probability, that is 
$$
 \lim_{m,n\to \infty} \PP\left( \A_{m\times n} \,\, \text{has} \,\,
   {\cal P}_{m,n} \right)=1
$$ 
is frequently considered, and -- by the Borel--Cantelli Lemma --
it implies almost sure convergence, if in addition 
$\sum_{m=1}^{\infty} \sum_{n=1}^{\infty} p_{mn} <\infty$ also holds,
where
$$p_{mn} = \PP\left( \A_{m\times n} \,\, \text{does not have} \,\,
   {\cal P}_{m,n} \right) .
$$


According to a generalization of a theorem of
F\"uredi and Koml\'os~\cite{FK} to rectangular matrices,
the spectral norm of an $m\times n$ Wigner-noise is $\sqrt{m+n}$ in
probability. More precisely, it was shown (see~\cite{AM}) that 
with probability tending to 1, 
$\| \W \| \le \frac73 \sigma  \sqrt{m+n}$, where $\sigma$ is the common
bound for the variances of the entries. Trivially, $\sigma \le K$ that
does not depend on $m$ and $n$, hence 
$\| \W \| = {\cal O} (\sqrt{m+n})$ in probability.
Bounding the variances from below, authors also proved that
$\| \W \| = \Theta (\sqrt{m+n})$ with high probability for large $m,n$.

To prove almost sure convergence, a 
 sharp concentration theorem of  N. Alon at al.
plays a crucial role (cf. \cite{AKV}). 
For completeness we formulate 
this result.



\begin{lem}\label{ev}
 Let $\widetilde \W$ be a $q\times q$ real symmetric matrix
whose entries in and above the main diagonal are independent
random variables with absolute value at most 1. Let $\lambda_1\ge
\lambda_2 \ge\dots \ge \lambda_q$ be the eigenvalues of 
$\widetilde \W$.
The following estimate holds for the deviation of the $i$th largest 
eigenvalue from its expectation with any positive real number $t$:
$$
 \PP \left( |\lambda_i - \EE (\lambda_i ) |>t \right) 
  \le \exp \left( -\frac{(1-o(1))t^2}{32i^2} \right)  \quad 
{\text{when}} \quad i\le \frac{q}2 ,
$$
and the same estimate holds for the probability
$ \PP \left( |\lambda_{q-i+1} - \EE (\lambda_{q-i+1} ) |>t \right) $.
\end{lem}

Now let $\W$ be a Wigner-noise with entries uniformly bounded by $K$. The 
\newline $(m+n)\times (m+n)$ symmetric matrix 
$$
 {\widetilde \W} =\frac1{K} \cdot
\begin{pmatrix} \0 & W \\
                 W^T & \0
\end{pmatrix}
$$
satisfies the conditions of Lemma~\ref{ev},  
its  largest  and  smallest  eigenvalues are 
$$
 \lambda_i ({\widetilde \W}) =-\lambda_{n+m-i+1} ({\widetilde \W}) =
 \frac1{K} \cdot s_i (\W ) , \qquad i=1,\dots ,\min \{ m,n \},
$$
the others are zeros,
where $\lambda_i (.)$ and $s_i (.)$ denote the $i$th largest
eigenvalue and singular value of the matrix in the argument,
respectively (cf. \cite{Bh}). Therefore  
\begin{equation}\label{s_1}
 \PP \left( | s_1 (\W ) - \EE (s_1 (\W ))|  > t \right) 
   \le \exp \left(- \frac{(1-o(1))t^2} {32 K^2} \right).    
\end{equation}
The fact that $\| \W \| = {\cal O} (\sqrt{m+n})$ in probability and 
inequality~(\ref{s_1}) together  ensure that 
$\EE(\Vert \W \Vert)={\cal O}(\sqrt{m+n})$. 
Hence, no matter how $\EE (\Vert \W \Vert)$ behaves
when $m\to \infty$
and $n\to \infty$, the following rough estimate holds.
\begin{lem}\label{lem-w-error}
There exist  positive constants $C_{K1}$
and $C_{K2}$,  
depending on the common bound on 
the entries of $\W$, such that
\begin{equation} \label{w-error}
\PP \left( \, \Vert \W\Vert > C_{K1} \cdot \sqrt{m+n}  \,  \right) 
   \le \exp[-C_{K2} \cdot (m+n)].   
\end{equation}
\end{lem}
The exponential decay of the right hand side of (\ref{w-error})
implies that the spectral norm of a Wigner-noise 
$\W_{m\times n}$   is of order $\sqrt{m+n}$, almost surely.
This observation will provide the base of almost sure
results of Sections~\ref{noisy} and~\ref{svpairs}.

In Section~\ref{noisy} we shall prove that the $m\times n$ noisy matrix $\A =\BB +\W$
almost surely has $r =\rk (\PPP )$ protruding singular values of
order $\sqrt{mn}$. In Section~\ref{svpairs} the distances of the corresponding
isotropic subspaces are estimated and this gives rise to a two-way
classification of the row and column items of $\A$ with sum of inner
variances ${\cal O} (\frac{m+n}{mn} )$, almost surely.  

In Definition~\ref{blownup} we required that the entries 
of the pattern matrix $\PPP$
be in the [0,1] interval. We made this restriction only for the sake of
the generalized Erd\H{o}s--R\'enyi hypergraph model 
with the entries of $\PPP$ as probabilities, see~\cite{bollobas}.  
In fact, our results are valid for any pattern matrix with fixed sizes 
and with non-negative entries. 
For example, in microarray measurements the rows 
correspond to different genes, the columns correspond to different conditions,
and the entries are the expression levels of a specific gene under a
specific condition.  

Sometimes the pattern matrix $\PPP$ is an $a\times b$ contingency table
with entries that are nonnegative integers.
Then the blown up matrix $\BB$
can be regarded as a larger ($m\times n$) contingency table that 
contains e.g., counts for two categorical variables with $m$ and
$n$ different categories, respectively. 
As the categories may be measured in different units, a normalization is
necessary. This normalization is made by dividing the entries of $\BB$
by the square roots of the corresponding row and column sums 
(cf.~\cite{KBCG}). 
This transformation is identical to that of the
correspondence analysis~\cite{G}, 
and the transformed matrix remains the same when we 
multiply  the initial matrix by a positive  constant. 
The transformed matrix $\BB_{corr}$, which belongs to $\BB$,  
has entries in [0,1] and maximum singular value 1. 
It is proved that there is a 
remarkable gap between the $\rk (\BB ) =\rk (\PPP )$ largest and the other
singular values of $\A_{corr}$, the matrix obtained from the  noisy
matrix $\A =\BB +\W$ by the correspondence transformation. 
This implies  well two-way classification properties 
of the row and column categories (genes and expression levels) 
in Section~\ref{corrm}. 

In Section~\ref{struct} a construction is given how a blown up structure behind a
real-life matrix with a few 
protruding singular values and 'well classifiable'
corresponding singular vector pairs can be found. 

\section{Singular values of a noisy matrix} \label{noisy}

\begin{prop}\label{blownup-sv}
If {\it GC1} holds, then all the non-zero singular values of the $m\times n$
 blown-up matrix $\BB$ are of order $\sqrt{mn}$.
\end{prop}

\begin{pf}
As there are at most $a$ and $b$ linearly
independent rows and linearly independent columns in $\BB$, respectively,
the rank $r$ of the matrix $\BB$ cannot exceed $\min \{ a,b \}$. 
Let $s_1 \ge s_2 \ge \dots \ge s_r >0$ be the positive singular values
of $\BB$.
Let $\vv_k \in \R^m$, $\uu_k \in \R^n$ be a singular vector pair
corresponding to
$s_k$, $k=1,\dots ,r$. 
Without loss of generality,  $\vv_1 ,\dots ,\vv_r$ and 
$\uu_1 ,\dots ,\uu_r$  can be  unit-norm, pairwise
orthogonal vectors in $\R^m$ and $\R^n$, respectively. 

For the subsequent calculations we drop the subscript  $k$, and $\vv$, $\uu$
denotes a singular vector pair 
corresponding to the singular value $s>0$
of the blown-up matrix $\BB$, $\| \vv \| =\| \uu \|=1$.
 It is easy to see that they have piecewise
constant structures: $\vv$
has $m_i$ coordinates equal to $v (i)$ $(i=1,\dots ,a)$ and $\uu$ has
$n_j$ coordinates equal to $u (j)$ $(j=1,\dots ,b)$. Then, with these
coordinates the singular value--singular vector equation 
\begin{equation}\label{Bu}
  \BB \uu = s \cdot \vv    
\end{equation}
has the form 
\begin{equation}\label{Bu-coord}
 \sum_{j=1}^b n_j p_{ij} u (j) =s \cdot v (i) \qquad (i=1,\dots ,a).   
\end{equation}
With the notations 
$$
 {\tilde \uu } =\left( u (1) ,\dots ,u (a) \right)^T , \qquad 
 {\tilde \vv }=\left( v (1) ,\dots ,v (b) \right)^T ,
$$
$$
 {\DD}_{m} = \diag (m_1 ,\dots ,m_a ), \qquad 
 {\DD}_{n} = \diag (n_1 ,\dots ,n_b ) 
$$ 
the equations in  (\ref{Bu-coord}) can be written as
$$
 \PPP {\DD}_{n} {\tilde \uu} =s \cdot {\tilde \vv}. 
$$ 
Introducing the following transformations of ${\tilde \uu}$ and ${\tilde \vv}$
\begin{equation}\label{w-z}
 \w = {\DD}^{1/2}_{n} {\tilde \uu} , \qquad
 \z = {\DD}^{1/2}_{m} {\tilde \vv} , 
\end{equation}
the equation is equivalent to 
\begin{equation}\label{sz}
 {\DD}_{m}^{1/2} {\PPP} {\DD}_{n}^{1/2} \w =s \cdot \z. 
\end{equation}
Applying the transformation~(\ref{w-z}) for the 
${\tilde \uu}_k , {\tilde \vv}_k$ pairs obtained from the 
$\uu_k ,\vv_k$ pairs $(k=1,\dots ,r)$,
orthogonormal systems in $\R^{a}$ and $\R^{b}$ are obtained:
$$
{\w_k}^T \cdot {\w_{\ell}}=\sum_{j=1}^b n_j u_k (j) u_{\ell} (j) =
\delta_{k\ell}
 \quad\text{and}\quad
  {\z_k}^T \cdot {\z_{\ell}}=\sum_{i=1}^a m_i v_k (i) v_{\ell} (i) =
\delta_{k\ell} .
$$
Consequently, $\z_k$, $\w_k$ is a singular vector pair corresponding
to  singular value $s_k$ of the $a\times b$ matrix ${\DD}_{m}^{1/2}
{\PPP} {\DD}_{n}^{1/2}$\quad ($k=1,\dots ,r$). With the shrinking 
$$
 {\widetilde \DD}_{m} = \frac1{m} {\DD}_{m}, \quad
 {\widetilde \DD}_{n} = \frac1{n} {\DD}_{n} 
$$ 
an equivalent form of (\ref{sz}) is 
$$
 {\widetilde \DD}_{m}^{1/2} {\PPP} {\widetilde \DD}_{n}^{1/2} \w =
    \frac{s}{\sqrt{mn}} \cdot \z , 
$$ 
that is the $a\times b$ matrix ${\widetilde
\DD}_{m}^{1/2} {\PPP} {\widetilde \DD}_{n}^{1/2}$ has non-zero singular
values $\frac{s_k}{\sqrt{mn}}$ with the same singular vector pairs $\z_k
,\w_k$\quad ($k=1,\dots ,r$).
If the $s_k$'s are not distinct numbers, the singular vector pairs
corresponding to a multiple singular value are not unique, but still
they can be obtained from the SVD  of the  shrunken matrix
${\widetilde \DD}_{m}^{1/2} {\PPP} {\widetilde \DD}_{n}^{1/2}$.

Now we want to establish relations between the singular values of $\PPP$
and ${\widetilde \DD}_{m}^{1/2} {\PPP} {\widetilde \DD}_{n}^{1/2}$.
Let $s_k(\QQ )$ denote the $k$th largest singular value of a matrix $\QQ$.
By the Courant--Fischer--Weyl minimax principle (cf.~\cite[p.75]{Bh})  
$$
  s_k(\QQ )=\max_{\dim H = k}\,\,\min_{\x \in H}\frac{\| \QQ \x \|}{\| \x \|}.
$$ 
Since we are interested only in the first $r$ singular values, where
$r = \rk(\BB) =
\rk({\widetilde \DD}_{m}^{1/2} {\PPP} {\widetilde \DD}_{n}^{1/2})$,
it is sufficient to consider vectors $\x$, for which
$ {\widetilde \DD}_{m}^{1/2} {\PPP} {\widetilde \DD}_{n}^{1/2}\x \ne \0$. 
Therefore with $k\in \{ 1,\dots ,r \}$ and an arbitrary $k$-dimensional
subspace $H\subset \R^b$ one can write
\begin{eqnarray*}
\min_{\x \in H} 
\frac{\| {\widetilde \DD}_{m}^{1/2} {\PPP} {\widetilde \DD}_{n}^{1/2} \x \| }{
\| \x \|} = 
\min_{\x \in H} 
\frac{\| {\widetilde \DD}_{m}^{1/2} {\PPP} {\widetilde \DD}_{n}^{1/2} \x \| }{
\| {\PPP} {\widetilde \DD}_{n}^{1/2} \x \| } \cdot 
\frac{\| {\PPP} {\widetilde \DD}_{n}^{1/2} \x \| }{
\| {\widetilde \DD}_{n}^{1/2} \x \|} \cdot
\frac{\| {\widetilde \DD}_{n}^{1/2} \x \|}{\| \x \|} \\
\ge s_a ({\widetilde \DD}_{m}^{1/2}) \cdot 
\min_{\x \in H}\frac{\| {\PPP} {\widetilde \DD}_{n}^{1/2} \x \| }{
\| {\widetilde \DD}_{n}^{1/2} \x \|} \cdot 
  s_b ({\widetilde \DD}_{n}^{1/2}) 
\ge \sqrt{cd} \cdot
\min_{\x \in H}\frac{\| {\PPP} {\widetilde \DD}_{n}^{1/2} \x \| }{
\| {\widetilde \DD}_{n}^{1/2} \x \|} ,
\end{eqnarray*}
with $c,d$ of {\it GC1}.
Now taking the maximum for all possible $k$-dimensional subspace $H$ we obtain
that $s_k({\widetilde \DD}_{m}^{1/2} {\PPP} {\widetilde \DD}_{n}^{1/2}) \ge 
\sqrt{cd} \cdot s_k(\PPP) >0$.
On the other hand, 
$$
s_k({\widetilde \DD}_{m}^{1/2} {\PPP} {\widetilde \DD}_{n}^{1/2}) \le
\| {\widetilde \DD}_{m}^{1/2} {\PPP} {\widetilde \DD}_{n}^{1/2} \| \le
\| {\widetilde \DD}_{m}^{1/2}  \| \cdot \| \PPP \| \cdot
\| {\widetilde \DD}_{n}^{1/2} \|
 \le \| \PPP \| \le \sqrt{ab} .
$$
These inequalities imply that 
$s_k({\widetilde \DD}_{m}^{1/2} {\PPP} {\widetilde \DD}_{n}^{1/2})$
is a nonzero constant, and because of 
$s_k({\widetilde \DD}_{m}^{1/2} {\PPP} {\widetilde \DD}_{n}^{1/2}) =
\frac{s_k}{\sqrt{mn}}$  we obtain that $s_1 ,\dots ,s_r =\Theta (\sqrt{mn})$.
\qed 
\end{pf}

\begin{thm}\label{svgap}
Let $\A =\BB +\W$ be an $m\times n$ random matrix,
where $\BB$ is a blown up matrix  with positive singular
values $s_1 ,\dots ,s_r$ and $\W$ is a Wigner-noise. 
Then, under {\it GC1}, the matrix $\A$ almost surely has $r$ singular values 
$z_1 , \dots ,z_r$,  such that 
$$
 |z_i -s_i| ={\cal O} (\sqrt{m+n}) , \qquad i=1,\dots ,r
$$ 
and for the other singular values almost surely
$$
  z_j ={\cal O} (\sqrt{m+n}) , \qquad j=r+1 ,\dots ,\min \{ m,n \}.
$$ 
\end{thm}

\begin{pf}
The statement  follows from the analog of the
Weyl's perturbation theorem for singular values of rectangular matrices
(see~\cite[p.99]{Bh}) and from Lemma~\ref{lem-w-error}.   
If $s_i (\A)$ and $s_i (\BB)$ denote the $i$th largest
singular values of the matrix in the argument 
then for the difference of the corresponding pairs
$$
 |s_i (\A ) -s_i (\BB )| \le \max_i s_i (\W ) =\| \W \| , \qquad
    i=1 ,\dots ,\min \{ m,n \}.
$$
By Lemma~\ref{lem-w-error},
$\PP \left( |s_i (\A ) -s_i (\BB )| > C_{K1} \cdot \sqrt{m+n}  \right) 
\le \PP \left( \, \Vert \W\Vert > C_{K1} \cdot \sqrt{m+n}  
\right) $
$   \le \exp[-C_{K2} \cdot (m+n)].$
The right hand side of the last inequality is  the general term of a 
convergent series (defined as a double summation),
thus  the convergence in probability 
implies the almost sure statement of the theorem. 
\qed
\end{pf} 

\begin{cor}
With notations 
\begin{equation}\label{eps-Delta}
 \eps :=\| \W \| ={\cal O} (\sqrt{m+n}) \quad \text{and} \quad
 \Delta := \min_{1\le i\le r} s_i (\BB ) =
 \min_{1\le i\le r} s_i = \Theta (\sqrt{mn})   
\end{equation} 
there is a spectral gap of size $\Delta -2\eps$ between the $r$
largest and the other singular values of the perturbed matrix $\A$,
and this gap is significantly larger than $\eps$. 
\end{cor}

\section{Classification via  singular vector pairs} 
\label{svpairs}

With the help of  Theorem~\ref{svgap}  we can  estimate the distances
between the corresponding right- and left-hand side eigenspaces 
(isotropic subspaces) of the  matrices $\BB$ and $\A =\BB+\W$.
Let $\vv_1 ,\dots ,\vv_m \in \R^m$ and $\uu_1 ,\dots ,\uu_n \in \R^n$ be
orthonormal left- and right-hand side singular vectors of $\BB$,
$$
 \BB \uu_i = s_i \cdot \vv_i \quad (i=1,\dots ,r) \quad \text{and} \quad
 \BB \uu_j = 0 \quad (j=r+1 ,\dots ,n).
$$
Let us also denote the unit-norm, pairwise orthogonal
left- and right-hand side singular vectors
corresponding to the $r$ protruding singular values $z_1 ,\dots ,z_r$ of
$\A$ by $\y_1 ,\dots ,\y_r \in \R^m$ and $\x_1 ,\dots ,\x_r \in \R^n$,
respectively. Then 
$\A \x_i = z_i \cdot \y_i$ \, $(i=1,\dots ,r)$. Let 
$$
 F:= \Span \{ \vv_1 , \dots ,\vv_r \}  \quad \text{and} \quad
 G:= \Span \{ \uu_1 , \dots ,\uu_r \} 
$$
denote the spanned linear subspaces in $\R^m$ and $\R^n$, respectively; 
further, let $\dist (\y ,F )$ denote the
Euclidean distance between the vector $\y$ and the subspace
$F$.

\begin{prop}\label{prop-dist}
 With the above notation, under {\it GC1}, 
the following estimate holds almost surely:
\begin{equation}\label{dist-F}
 \sum_{i=1}^r {\dist}^2 (\y_i ,F ) \le r \frac{\eps^2}{(\Delta -\eps )^2 }
 ={\cal O}\left(\frac{m+n}{mn}\right)
\end{equation}
and analogously, 
\begin{equation}\label{dist-G}
 \sum_{i=1}^r {\dist}^2 (\x_i ,G ) \le r \frac{\eps^2}{(\Delta -\eps )^2 }
 ={\cal O}\left(\frac{m+n}{mn}\right).   
\end{equation}
\end{prop}

\begin{pf}
Let us choose one of  the right-hand side singular vectors $\x_1 ,\dots ,\x_r$
of $\A =\BB +\W$ and denote it simply by $\x$ with corresponding
singular value  $z$. We shall estimate the
distance between $\x$ and $G$, similarly between $\y =\A \x /z$ and $F$. 
For this purpose we expand $\x$ and $\y$ in the orthonormal
bases $\uu_1 ,\dots ,\uu_n$ and $\vv_1 ,\dots ,\vv_m$, respectively:
$$
  \x =\sum_{i=1}^n t_i \uu_i \quad \text{and} \quad
  \y =\sum_{i=1}^m l_i \vv_i .
$$
Then
\begin{equation}\label{Ax-1}
  \A \x = (\BB +\W )\x = \sum_{i=1}^r t_i s_i \vv_i  + \W \x , 
\end{equation}
and, on the other hand,
\begin{equation}\label{Ax-2}
 \A \x = z \y = \sum_{i=1}^m  z l_i \vv_i . 
\end{equation}
Equating the right-hand sides of (\ref{Ax-1}) and (\ref{Ax-2})  we obtain
$$
  \sum_{i=1}^r (z l_i -t_i s_i ) \vv_i  
 +\sum_{i=r+1}^{m} z l_i \vv_i = \W \x .
$$
Applying the Pythagorean Theorem
\begin{equation}
 \sum_{i=1}^r (z l_i -t_i s_i )^2 
 + z^2 \sum_{i=r+1}^m l_i^2   
 = \|\W \x\|^2 \le   \eps^2 ,    
\end{equation}
because $\| \x \| =1$ and  $\|\W \| =  \eps$.

As $z \ge \Delta -\eps$  holds almost surely by Theorem~\ref{svgap},
$$
  {\dist}^2 (\y ,F)= \sum_{i=r+1}^m l_i^2 \le \frac{\eps^2}{z^2}
  \le \frac{\eps^2}{(\Delta-\eps)^2}.
$$
The order of the above estimate follows from the order of $\eps$ and 
$\Delta$ of (\ref{eps-Delta}):
\begin{equation}\label{dist-F-eps}
 {\dist}^2 (\y ,F) 
= {\cal O}(\frac{m+n}{mn})          
\end{equation}
almost surely. Applying (\ref{dist-F-eps}) for the left-hand side singular 
vectors $\y_1 , \dots ,\y_r$, by the Definition~\ref{conv} 
\begin{multline*}
\PP\left\{
  \exists  m_{0i}, n_{0i} \in \N \,\, \text{such that for} 
\,\,
 m\ge m_{0i} \,\, \text{and} \,\, n\ge n_{0i} \colon \right.
\\
 \left. {\dist}^2 (\y_i ,F) \le  {\eps^2 }/{(\Delta -\eps )^2} 
\right\} =1
\end{multline*}
for $i=1,\dots ,r$. Hence,
\begin{multline*}
\PP\left\{ \exists  m_0, n_0 \in \N \,\, \text{such that for} \,\,
 m\ge m_0 \,\, \text{and} \,\, n\ge n_0 \colon \right. \\
 \left. {\dist}^2 (\y_i ,F) \le {\eps^2 }/{(\Delta -\eps )^2}, \, 
i=1,\dots , r  \right\} =1 ,
\end{multline*}
consequently,
\begin{multline*}
\PP\left\{ \exists  m_0, n_0 \in \N \,\, \text{such that for} \,\,
 m\ge m_0 \,\, \text{and} \,\, n\ge n_0  \colon \right. \\
  \sum_{i=1}^r {\dist}^2 (\y_i ,F) \le r {\eps^2 }/{(\Delta 
-\eps )^2}\,\} =1 
\end{multline*}
also holds, and this finishes the proof of the first statement.

The estimate for the squared distance between $G$ and a right-hand side 
singular vector $\x$ of $\A$ follows in the same way starting with
$\A^T \y = z \cdot \x $
and using the fact that $\A^T$ has the same singular values as $\A$.
\qed
\end{pf}

By Proposition~\ref{prop-dist},  the individual distances 
between the original and the perturbed subspaces and also
the sum of these distances tend to zero almost surely as $m,n\to\infty$.

Now let $\A$ be a microarray on  $m$ genes and $n$ 
conditions, with $a_{ij}$ denoting the expression level of gene $i$ under
condition $j$.  
We suppose that $\A$  is a noisy
random matrix obtained by adding a Wigner-noise $\W$ to the blown up
matrix $\BB$. Let us denote by $A_1 ,\dots ,A_a$ the partition of the genes
and by $B_1 ,\dots ,B_b$ the partition of the conditions
with respect to the blow-up (they can also be thought of as clusters of
genes and conditions).

Proposition~\ref{prop-dist}
also implies the well-clustering property of the representatives
of the genes and conditions in the following representation. Let $\Y$ be the
$m\times r$ matrix containing the left-hand side singular vectors 
$\y_1 ,\dots ,\y_r$ of $\A$ in its columns. 
Similarly, let $\X$ be the
$n\times r$ matrix containing the right-hand side singular vectors 
$\x_1 ,\dots ,\x_r$ of $\A$ in its columns. 
Let the $r$-dimensional representatives of the genes be
the row vectors of $\Y$: \, $\y^{1}, \dots ,\y^{m} \in \R^r$,
while the $r$-dimensional representatives of the conditions be
the row vectors of $\X$: \, $\x^{1}, \dots ,\x^{n} \in \R^r$.
 Let $S_a^2 (\Y )$ denote  the $a$-variance, introduced in 
\cite{bolla1}, 
of the genes' representatives 
$$
 S_a^2 (\Y )= \min_{ \{ A'_1 ,\dots ,A'_a \} } 
 \sum_{i=1}^a \sum_{j\in A'_i} \|\y^{j}- {\bar \y}^{i}\|^2 ,
 \quad \text{where} \quad
 {\bar \y}^{i} =\frac1{m_i} \sum_{j\in A'_i} \y^{j} ,  
$$
while $S_b^2 (\X )$ denotes  the $b$-variance 
of the conditions' representatives 
$$
 S_b^2 (\X )= \min_{ \{ B'_1 ,\dots , B'_b \} }
 \sum_{i=1}^b \sum_{j\in B'_i} \|\x^{j}- {\bar \x}^{i}\|^2 ,
 \quad \text{where} \quad
 {\bar \x}^{i} =\frac1{n_i} \sum_{j\in B'_i} \x^{j} ,  
$$
the partitions $\{ A'_1 ,\dots ,A'_a \}$ and 
$ \{ B'_1 ,\dots , B'_b \} $ varying over all 
$a$- and $b$-partitions of the genes and conditions, respectively.

\begin{thm}\label{thm-SA}
With the above notation, under {\it GC1},
for the $a$- and $b$-variances 
of the representation of the microarray $\A $  the relations
$$  
  S_a^2 (\Y ) = {\cal O} \left( \frac{m+n}{mn} \right) \quad \text{and} \quad 
  S_b^2 (\X ) = {\cal O} \left( \frac{m+n}{mn} \right)
$$
hold almost surely.
\end{thm}

\begin{pf}
 By the proof of Theorem 3 of~\cite{bolla1}
 it can be easily seen that 
$S_a^2 (\Y ) \le  \sum_{i=1}^a \sum_{j\in A_i} \|\y^{j}- {\bar \y}^{i}\|^2$
and 
$S_b^2 (\X ) \le \sum_{i=1}^b \sum_{j\in B_i} \|\x^{j}- {\bar \x}^{i}\|^2 ,$ 
the right-hand sides being equal to the 
left-hand sides of (\ref{dist-F}) and (\ref{dist-G}), respectively, 
therefore they are also of order $\frac{m+n}{mn}$. 
\qed
\end{pf}

Hence, the addition of  any kind of a Wigner-noise to  a rectangular matrix 
that has a blown up structure $\BB$ will not change 
the order of the protruding singular values, 
and the block structure of $\BB$ can be reconstructed
from the representatives of the row and column items of the noisy matrix $\A$.

With an appropriate Wigner-noise,  we
can achieve that the matrix $\BB +\W$  in its  $(i,j)$-th block 
contains 1's with probability $p_{ij}$, and 0's otherwise. That is,
for $i=1,\dots ,a$, \, $j=1,\dots ,b$,\,\, $l\in A_i$,\,\, $k\in B_j$,
 let
\begin{equation}\label{01-matrix}
 w_{lk} := \left\{
      \begin{array}{ll}
   1-p_{ij}, & \mbox{with probability } \quad p_{ij}  \\
   -p_{ij}   & \mbox{with probability } \quad 1-p_{ij} 
    \end{array}
\right.   
\end{equation}
 be  independent  random variables.
This $\W$ satisfies the conditions of Definition~\ref{wigner} with  entries 
uniformly bounded by 1, zero expectation and variance 
$$
 \sigma^2 =\max_{1\le i\le a; \, 1\le j\le b} p_{ij} (1-p_{ij}) 
\le \frac{1}{4}.
$$   
The noisy matrix $\A$ becomes a 0-1 matrix that can be regarded as the 
incidence matrix
of a hypergraph on $m$ vertices and $n$ edges. 
(Vertices correspond to the genes and edges correspond to the conditions.
The incidence relation depends on  
whether a specific gene is expressed 
or not under a specific condition).

By the choice (\ref{01-matrix}) of $\W$,
 vertices of the vertex set $A_i$ appear in edges of the  edge set $B_j$ 
with probability $p_{ij}$ (set $i$ of genes equally influences set $j$
of conditions, like the chess-board pattern of~\cite{KBCG}).
It is  a generalization of the classical Erd\H os--R\'enyi
model for random hypergraphs   and for several blocks,  
see~\cite{bollobas}.
The question, how such a chess-board pattern behind a random (especially
0-1) matrix can be found under specific conditions, is discussed in
Section~\ref{struct}.

\section{Perturbation results for correspondence matrices} \label{corrm}

Now the pattern matrix $\PPP$ contains 
arbitrary non-negative entries, so does the
blown up matrix $\BB$. Let us suppose that there are no identically zero
rows or columns. We perform the correspondence transformation described below
on $\BB$.
We are interested in the order of singular values of 
matrix $\A = \BB+\W$ when the same correspondence transformation is 
applied to it.
To this end, we introduce the following notations: 
\begin{eqnarray*}
& {\DD}_{Brow}&=
\diag(d_{Brow \, 1},\dots, d_{Brow \, m}):=
\diag 
\left(\sum_{j=1}^n b_{1j},\dots , \sum_{j=1}^n  b_{mj}\right)\\
&{\DD}_{Bcol}&=\diag(d_{Bcol \, 1},\dots, d_{Bcol \, n}):=
 \diag 
\left(\sum_{i=1}^m b_{i1}, \dots , \sum_{i=1}^m  b_{in}\right)
\\
& {\DD}_{Arow}&=\diag(d_{Arow \, 1},\dots, d_{Arow \, m}):= 
 \diag \left(\sum_{j=1}^n a_{1j}, \dots ,\sum_{j=1}^n 
 a_{mj}\right)
\\
& {\DD}_{Acol}&= \diag(d_{Acol \, 1},\dots, d_{Acol \, n}):=
\diag \left(\sum_{i=1}^m a_{i1}, \dots , \sum_{i=1}^m 
 a_{in}\right).
\end{eqnarray*}
Further, set
$$
\BB_{corr}:={\DD}_{Brow}^{-1/2}\BB {\DD}_{Bcol}^{-1/2} \quad\text{and}\quad
\A_{corr}:={\DD}_{Arow}^{-1/2}\A {\DD}_{Acol}^{-1/2}  
$$ 
for the transformed matrices obtained from $\BB$ and $\A$ while
carrying out correspondence analysis on $\BB$ and the same correspondence
transformation on $\A$.
It is well known~\cite{G} that the leading singular value  of $\BB_{corr}$ 
is equal to  1 and the multiplicity of 1 as a singular value coincides
with the number of irreducible blocks in $\BB$. Let $s_i$ denote a non-zero
singular value of $\BB_{corr}$ with unit-norm singular vector pair 
$\vv_i$, $\uu_i$. With the transformations
$$
  \vv_{corr \, i} :={\DD}^{-1/2}_{Brow} \vv_i \quad \text{and} \quad 
  \uu_{corr \, i} :={\DD}^{-1/2}_{Bcol} \uu_i  
$$
the so-called correspondence vector pairs 
are obtained.
If the coordinates $v_{corr \, i} (j)$, $u_{corr \, i} (\ell )$  
of such a pair 
are regarded as  possible values of
two discrete random variables $\beta_i$ and $\alpha_i$ (often called the 
$i$th correspondence factor pair)
with the prescribed marginals, then, as
in canonical analysis, their correlation is $s_i$, and this is the largest
possible correlation under the condition that they are uncorrelated with
the previous random variables $\beta_1 , \dots ,\beta_{i-1}$ and 
$\alpha_1 , \dots ,\alpha_{i-1}$,
respectively ($i>1$).

If $s_1 =1$ is a single singular value,
then $\vv_{corr \, 1}$ and $\uu_{corr \, 1}$ are the  all 1 vectors 
and the corresponding $\beta_1$, $\alpha_1$ pair is regarded as a
trivial correspondence factor pair. This corresponds to the
general case. 
Keeping $k\le \rk(\BB_{corr}) =\rk(\BB)  =\rk(\PPP)$ singular values with the
coordinates of the corresponding $k-1$ non-trivial correspondence factor
pairs, the following  $(k-1)$-dimensional representation of the 
$j$th and $\ell$th categories
of the underlying two discrete variables is obtained:
$$
 \vv^{j}_{corr} := \left( v_{corr \,2} (j) , \dots , v_{corr \, k} (j) \right) 
                      \quad \text{and} \quad 
 \uu^{\ell}_{corr} := \left( u_{corr \,2} (\ell ) , \dots , 
    u_{corr \, k} (\ell ) \right).
$$

This representation
has the following optimality properties: the closeness of categories of 
the same variable reflects the  similarity between them, while the closeness
of  categories of different variables reflects their frequent 
simultaneous occurrence.
For example, $\BB$ being a microarray, the representatives of similar
function genes, as well as representatives of similar conditions are close
to each other; also, representatives of genes that are responsible for a given
condition, are close to the representatives of those conditions.
Now we prove the following.

\begin{prop}\label{Bcorr-sv}
 Given the blown up matrix $\BB$, under {\it GC1}
 there exists a  constant $\delta \in (0,1)$, independent 
of $m$ and $n$, 
 such that  all the $r$ non-zero singular values
of  $\BB_{corr}$ are in the interval $[\delta , 1]$, where 
$r=\rk(\BB) =\rk(\PPP)$.
\end{prop}

\begin{pf}
It is easy to see that $\BB_{corr}$ is the blown up matrix of the $a\times b$
pattern matrix $\tilde{\PPP}$ with entries
$$
 \tilde{p}_{ij} = \frac{p_{ij}}
 {\sqrt{(\sum_{\ell =1}^b p_{i\ell} n_{\ell} ) (\sum_{k=1}^a p_{kj} m_k )}} .
$$  
Following the considerations of the proof of 
Proposition~\ref{blownup-sv},
the blown up matrix $\BB_{corr}$ 
has exactly $r=\rk(\PPP) =\rk(\tilde{\PPP})$ non-zero singular values 
that are the singular values of the $a\times b$ matrix
$\PPP' =\DD_m^{1/2} \tilde{\PPP} \DD_n^{1/2}$ with entries
$$
 p'_{ij} = 
 \frac{p_{ij} \sqrt{m_i} \sqrt{n_j} }
  {\sqrt{(\sum_{\ell =1}^b p_{i\ell} n_{\ell} ) 
(\sum_{k=1}^a p_{kj} m_k )}} =
 \frac{p_{ij}}{\sqrt{(\sum_{\ell =1}^b p_{i\ell} \frac{n_{\ell}}{n_j} ) 
   (\sum_{k=1}^a p_{kj} \frac{m_k}{m_i} )}} .
$$
Since the matrix $\PPP$ contains no identically zero rows or columns,
the matrix  $\PPP'$ varies on a compact set of $a\times b$ matrices
 determined by the inequalities (\ref{const}). 
The range of the non-zero
singular values depends continuously on the matrix that does not depend
on $m$ and $n$. Therefore, the minimum
non-zero singular value does not depend  on $m$ or $n$.
Because the largest singular value is 1, this finishes the proof. 
\qed
\end{pf}

\begin{thm} \label{Acorr-sv}
 Under {\it GC1} and {\it GC2}, there exists a positive 
number $\delta$ (independent
of $m$ and $n$) such that  for every $0< \tau <1/2$
the following statement holds almost surely: the  $r$ largest 
singular values of $\A_{corr}$
are in the interval $[\delta -\max \{ n^{-\tau }, m^{-\tau } \} , 1+ 
\max \{ n^{-\tau }, m^{-\tau } \} ]$, while all  the others 
are at most $\max \{ n^{-\tau }, m^{-\tau } \}$.
\end{thm}

\begin{pf}
 First notice that 
\begin{equation}\label{A-corr}
\A_{corr}={\DD}_{Arow}^{-1/2}\A {\DD}_{Acol}^{-1/2}=
 {\DD}_{Arow}^{-1/2}\BB {\DD}_{Acol}^{-1/2}+  
 {\DD}_{Arow}^{-1/2}\W {\DD}_{Acol}^{-1/2}.   
\end{equation}
The entries of ${\DD}_{Brow}$ and those of
${\DD}_{Bcol}$ are  of order $\Theta(n)$ and $\Theta(m)$, respectively.
Now we prove that  for every $i=1,\dots, m$ and
 $j=1,\dots, n$
$\vert d_{Arow\,i} -d_{Brow\,i} \vert < n\cdot n^{-\tau } $  and
$\vert d_{Acol\,j} -d_{Bcol\,j} \vert < m\cdot  m^{-\tau } $ hold 
almost surely.  To this end, we use Chernoff's
inequality for large deviations (cf.~\cite{bolla2}, Lemma 4.2):
\begin{equation*}
\begin{split}
\PP&\left(\vert d_{Arow\,i} -d_{Brow\,i} \vert > n\cdot n^{-\tau }
\right)  =
\PP\left( \left\vert \sum_{j=1}^n w_{ij} \right\vert > n^{1-\tau }
\right)  \\
&< \exp
\left\{-\frac{n^{2-2\tau }}{2 (\Var (\sum_{j=1}^n w_{ij})
 +Kn^{1-\tau }/3 )}\right\} 
\le \exp
\left\{-\frac{n^{2-2\tau } }{2 (n\sigma^2 +Kn^{1-\tau } /3 
)}\right\}\\
& =
 \exp\left\{-\frac{n^{1-2\tau } }{2 (\sigma^2 +K n^{-\tau } /3 
)}\right\} 
\quad (i=1,\dots ,m),
\end{split}
\end{equation*}
where the constant $K$ is the uniform bound for $|w_{ij}|$'s and $\sigma^2$
is the bound for their variances.
In virtue of {\it GC2} the following estimate holds 
with some $C_0 >0$ and $C\ge 1$ (constants of {\it GC2}) and
large enough  $n$:
\begin{equation}\label{d_row}
\begin{split}
\PP&\left( \vert d_{Arow\,i} -d_{Brow\,i} \vert >n^{1-\tau } \text{   
for all}
\quad i\in \{ 1,\dots ,m \} \right)  \\
 & \le 
m\cdot
 \exp\left\{-\frac{n^{1-2\tau } }{2 (\sigma^2 +K n^{-\tau } /3 
   )}\right\}   
  \le C_0 \cdot n^C \cdot
 \exp\left\{-\frac{n^{1-2\tau } }{2 (\sigma^2 +K n^{-\tau } /3 
   )}\right\}   \\
 & =  \exp\left\{ \ln C_0 + C\ln n
-\frac{n^{1-2\tau } }{2 (\sigma^2 +K n^{-\tau} /3  )}\right\} . 
\end{split}
\end{equation}
The estimation of probability 
$$
 \PP \left( \vert d_{Acol\,j} -d_{Bcol\,j} \vert > m^{1-\tau } \text{  
 for all} \quad j\in \{ 1,\dots ,n \} \right) 
$$
can be treated analogously (with  $D_0 >0$ and $D\ge 1$ of {\it GC2}). 
The right-hand side
of (\ref{d_row}) forms a convergent series, therefore 
\begin{equation}\label{min-d}
 \min_{i\in\{1,\dots, m\}}| d_{Arow\,i} | =\Theta (n), \qquad
  \min_{j\in\{1,\dots, n\}}| d_{Acol\,j} | =\Theta (m)
\end{equation}
hold  almost surely.

Now it is straightforward to bound the norm of the second 
term of (\ref{A-corr}) by
\begin{equation}\label{DWD}
 \| {\DD}_{Arow}^{-1/2} \| \cdot \|\W\| \cdot \| {\DD}_{Acol}^{-1/2} \| .
\end{equation}
As by Lemma~\ref{lem-w-error},  
$\| \W \| ={\cal O}(\sqrt{m+n})$ holds almost surely, the quantity (\ref{DWD}) is at 
most of order $\sqrt{\frac{m+n}{mn}}$ almost surely.
Hence, it is almost surely less than  $\max \{ n^{-\tau }, 
m^{-\tau } \}$.

To estimate the norm of the first term of (\ref{A-corr}) 
let us write it in the form
\begin{equation}\label{DBD}
\begin{split}
 {\DD}_{Arow}^{-1/2} \BB {\DD}_{Acol}^{-1/2}=  
 {\DD}_{Brow}^{-1/2} \BB {\DD}_{Bcol}^{-1/2} 
 &+\left[{\DD}_{Arow}^{-1/2}-{\DD}_{Brow}^{-1/2} \right] \BB
 {\DD}_{Bcol}^{-1/2} \\  
 &+ {\DD}_{Arow}^{-1/2} \BB
\left[{\DD}_{Acol}^{-1/2}-{\DD}_{Bcol}^{-1/2} \right] .
\end{split}
\end{equation}

 The first term  is just $\BB_{corr}$, so due to
Proposition~\ref{Bcorr-sv}, we should prove only that the norms 
of both remainder terms are almost surely less than 
$\max \{ n^{-\tau}, m^{-\tau} \}$.
These two terms have a similar appearance,
therefore it is enough to  estimate
one of them. For example, the second term can be bounded by 
\begin{equation}\label{normbound}
 \| {\DD}_{Arow}^{-1/2}-{\DD}_{Brow}^{-1/2} \| \cdot \| \BB \| \cdot
 \| {\DD}_{Bcol}^{-1/2} \| .   
\end{equation}
The  estimation of the first factor in (\ref{normbound}) is 
as follows:
\begin{equation}\label{Derror}
\begin{split}
 & \| {\DD}_{Arow}^{-1/2}-{\DD}_{Brow}^{-1/2} \| =\max_{i\in \{ 1,\dots 
,m \} }
\left( \frac1{\sqrt{d_{Arow\,i}}} - \frac1{\sqrt{d_{Brow\,i}}} \right) 
\\
&=\max_{i\in \{ 1,\dots ,m \} }
\frac{ |d_{Arow\,i} - d_{Brow\,i} |} 
{ \sqrt{ d_{Arow\,i} \cdot d_{Brow\,i} } 
(\sqrt{d_{Arow\,i}}+\sqrt{d_{Brow\,i}}) }   \\
&\le \max_{i\in \{ 1,\dots ,m \} }
\frac{ |d_{Arow\,i} - d_{Brow\,i} |} 
{ \sqrt{ d_{Arow\,i} \cdot d_{Brow\,i} }  } \cdot
  \max_{i\in \{ 1,\dots ,m \} }
\frac{ 1 } 
{  (\sqrt{d_{Arow\,i}}+\sqrt{d_{Brow\,i}}) } .
\end{split}
\end{equation}  
By relations (\ref{min-d}),
 $\sqrt{ d_{Arow\,i} \cdot d_{Browi} } =\Theta (n)$
for any $i=1,\dots ,m$, and hence,
$$
 \frac{ |d_{Arow\,i} - d_{Brow\,i} |}{\sqrt{ d_{Arow\,i} \cdot 
d_{Brow\,i} } } \le n^{-\tau}
$$
almost surely,   further
$\max_{i\in\{1,\dots, m\}}
\frac1{\sqrt{d_{Arowi}}+\sqrt{d_{Browi}}} =\Theta (\frac1{\sqrt{n}} )$
almost surely. 

Therefore the left hand side of (\ref{Derror}) can be estimated by 
$n^{-\tau -1/2}$ from above 
almost surely. 
For the further factors in (\ref{normbound}) we obtain
$\| \BB \| =\Theta ( \sqrt{mn } ) $ 
(see Proposition~\ref{blownup-sv}), while
$\| {\DD}_{Bcol}^{-1/2} \| =\Theta (\frac1{\sqrt{m}})$ almost surely.
These together imply that
 $$
 n^{-\tau -1/2}\cdot n^{1/2}m^{1/2}\cdot m^{-1/2} \le n^{-\tau} 
 \le \max \{ n^{-\tau}, m^{-\tau} \}. 
$$
This finishes the estimation of the first term in 
(\ref{A-corr}), and 
by he Weyl's perturbation theorem the proof, too.  
\qed
\end{pf}

\begin{rem}
In the Gaussian case  the large deviation principle 
can be replaced by the  simple estimation of the Gaussian probabilities
with any $\kappa >0$:
$$
\PP\left( \left\vert\frac 1n \sum_{j=1}^n w_{ij} \right\vert> \kappa
\right)< \min\left(1, \frac {4\sigma} {\kappa \sqrt{2\pi n}}\exp\left\{
- \frac{n} {2\sigma^2} \kappa^2 \right\}\right).  
$$
Setting $\kappa =n^{-\tau}$ we get an estimate, analogous to 
(\ref{d_row}).
\end{rem}

Suppose that the blown up matrix $\BB$ is irreducible and its non-negative
entries sum up to 1. This restriction does not effect the result of the
correspondence analysis, that is the SVD of the matrix $\BB_{corr}$.
Remember that the non-zero singular values of
$\BB_{corr}$ are the numbers $1=s_1 >s_2 \ge \dots \ge s_r >0$ 
with unit-norm singular vector pairs 
$\vv_i$, $\uu_i$ having piecewise constant structure ($i=1,\dots ,r)$.
Set
$$ 
 F:= \Span \{ \vv_1 , \dots ,\vv_r \}  \quad \text{and} \quad 
 G:= \Span \{ \uu_1 , \dots ,\uu_r \} .
$$ 

Let $0<\tau <1/2$ be arbitrary and
$\epsilon :=\max \{ n^{-\tau }, m^{-\tau } \}$.  
Let us also denote the unit-norm, pairwise orthogonal 
left- and right-hand side singular vectors
corresponding to the $r$ singular values $z_1 ,\dots ,z_r \in 
[\delta -\epsilon ,1+\epsilon ]$ of 
$\A_{corr}$ -- guaranteed by Theorem~\ref{Acorr-sv} 
under {\it GC2}  -- 
by $\y_1 ,\dots ,\y_r \in \R^m$ and $\x_1 ,\dots ,\x_r \in \R^n$,
respectively.
 
\begin{prop}\label{dist-FG-est}
With the above notation, under {\it GC1} and
{\it GC2} the following estimate holds almost surely 
for the  distance between $\y_i$ and $F$:
\begin{equation}\label{dist-F-est}
  {\dist} (\y_i ,F ) \le  
 \frac{\epsilon}{(\delta -\epsilon ) }
 =\frac1{(\frac{\delta}{\epsilon} -1)}
 \qquad (i=1,\dots ,r)   
\end{equation}
and analogously, for the distance between  $\x_i$ and $G$:
\begin{equation}\label{dist-G-est}
  {\dist} (\x_i ,G ) \le 
 \frac{\epsilon}{(\delta -\epsilon )}
 =\frac1{(\frac{\delta}{\epsilon} -1)} 
 \qquad (i=1,\dots ,r) .   
\end{equation}
\end{prop}

\begin{pf}
Follow the method of proving Proposition~\ref{prop-dist} 
-- under {\it GC1} -- with 
$\delta$ instead of $\Delta$ and $\epsilon$ instead of $\eps$. 
Here {\it GC2} is necessary only for $\A_{corr}$ to have $r$ protruding
singular values.
\qed
\end{pf}

\begin{rem}
The left-hand sides of (\ref{dist-F-est}) and (\ref{dist-G-est}) are 
almost surely of order $\max \{ n^{-\tau}, m^{-\tau} \}$
that tend to zero as  $m,n\to \infty$ under {\it GC1} and {\it GC2}.
\end{rem}

Proposition~\ref{dist-FG-est} 
implies the well-clustering property of the representatives
of the two discrete variables 
by means of the noisy correspondence vector pairs
$$
  \y_{corr \, i} :={\DD}^{-1/2}_{Arow} \y_i , \qquad 
  \x_{corr \, i} :={\DD}^{-1/2}_{Acol} \x_i \quad (i=1,\dots ,r).
$$
Let $\Y_{corr}$  denote the
$m\times r$ matrix that contains the left-hand side vectors 
$\y_{corr \, 1} ,\dots ,\y_{corr \, r}$ in its columns. 
Similarly, let $\X_{corr}$ denote the
$n\times r$ matrix that contains the right-hand side vectors
$\x_{corr \, 1} ,\dots ,\x_{corr \, r}$ in its columns. 
The $r$-dimensional representatives of $\alpha$  are
the row vectors of $\Y_{corr}$ denoted by
  $\y^{1}_{corr}, \dots ,\y^{m}_{corr} 
\in \R^r$, while the $r$-dimensional representatives of $\beta$ are
the row vectors of $\X_{corr}$ denoted by $\x^{1}_{corr}, \dots ,\x^{n}_{corr} 
\in \R^r$. With respect to the marginal
distributions, let the  $a$- and $b$-variances of these 
representatives be defined by
$$
 S_a^2 (\Y_{corr} )= \min_{ \{ A'_1 ,\dots ,A'_a \} }
\sum_{i=1}^a \sum_{j\in A'_i} d_{Arow \, j} 
 \|\y^{j}_{corr}- {\bar \y}^{i}_{corr} \|^2 \ ,$$
$$
 S_b^2 (\X_{corr} )= \min_{ \{ B'_1 ,\dots , B'_b \} }
\sum_{i=1}^b \sum_{j\in B'_i} d_{Acol \, j} 
 \|\x^{(j)}_{corr}- {\bar \x}^{i}_{corr} \|^2 ,
$$
where $\{ A'_1 ,\dots ,A'_a \}$ and $ \{ B'_1 ,\dots , B'_b \} $ 
are
$a$- and $b$-partitions of the genes and conditions, respectively, 
$$ {\bar \y}^{i}_{corr} =\sum_{j\in A'_i} d_{Arow \, j} 
\y^{j}_{corr} \quad\text{ and }\quad 
 {\bar \x}^{i}_{corr} =\sum_{j\in B'_i} d_{Acol \, j} \x^{j}_{corr}.
$$ 

\begin{thm} \label{thm-S}
 With the above notation, under {\it GC1} and
{\it GC2},
$$
 S_a^2 (\Y_{corr} ) \le \frac{r}{(\frac{\delta}{\epsilon} -1)^2}  
    \qquad \text{and} \qquad
 S_b^2 (\X_{corr} ) \le \frac{r}{(\frac{\delta}{\epsilon} -1)^2}  
$$
hold  almost surely, where
$\epsilon =\max \{ n^{-\tau }, m^{-\tau } \}$ with 
every $0<\tau <1/2$.
\end{thm}

\begin{pf}
An easy calculation  shows  that 
$$
S_a^2 (\Y_{corr} ) \le
\sum_{i=1}^a \sum_{j\in A_i} d_{Arow \, j} 
 \|\y^{j}_{corr}- {\bar \y}^{i}_{corr} \|^2 =
 \sum_{i=1}^{r} {\dist}^2 (\y_i ,F ), 
$$
$$
S_b^2 (\X_{corr} ) \le
\sum_{i=1}^b \sum_{j\in B_i} d_{Acol \, j} 
 \|\x^{(j)}_{corr}- {\bar \x}^{i}_{corr} \|^2 =
 \sum_{i=1}^{r} {\dist}^2 (\x_i ,G ) , 
$$ 
hence the result of Proposition~\ref{dist-FG-est} can be used.
\qed
\end{pf}

Under {\it GC1} and
{\it GC2} with $m,n$ large enough, Theorem~\ref{thm-S} implies that after
 performing correspondence analysis on the noisy matrix
$\A$, the representation through the correspondence vectors
belonging to $\A_{corr}$ will also reveal the block structure behind $\A$.

\section{Recognizing the structure} \label{struct}

One might wonder where the singular values of an
$m\times n$ matrix $\A =(a_{ij})$ are located  if
 $a:=\max_{i,j} | a_{ij}|$ is independent of $m$ and $n$.
On one hand, the maximum singular value 
cannot exceed ${\cal O} (\sqrt{mn})$, as it is
at most $\sqrt{ \sum_{i=1}^m \sum_{j=1}^n a_{ij}^2 }$.
On the other hand, let
$\QQ$ be an $m\times n$ random matrix with entries $a$ or $-a$
(independently of each other).
Consider the spectral norm of all such matrices and take the minimum
of them:
$
  \min_{\QQ \in \{ -a, +a \}^{m\times n} } \| \QQ \| .
$
\ This quantity measures the minimum linear structure that a matrix  of the
same size and magnitude as $\A$ can possess. As the Frobenius
norm  of  $\QQ$ is $a \sqrt{mn} $,
in virtue of inequalities between spectral and Frobenius norms,
the above minimum is at least $\frac{a}{\sqrt{2}} \sqrt{m+n}$,
which is exactly the order of the spectral norm of a Wigner-noise.

So an $m\times n$ random matrix (whose entries are independent
and uniformly bounded)
under very general conditions has at least one singular value of order
greater than $\sqrt{m+n}$.
Suppose there are $k$ such singular values
and the representatives by means of the corresponding singular vector
pairs can be well classified in the sense of Theorem~\ref{thm-SA}
(cf. the introduction to that theorem).  
 Under these conditions we can reconstruct
a blown up structure behind our matrix.

\begin{thm}\label{constr}
Let  $\A_{m\times n}$  be a  sequence of $m\times n$ matrices, where $m$ 
and $n$ tend to infinity.
 Assume, that $\A_{m\times n}$ has exactly $k$ 
singular values 
of order  greater than $\sqrt{m+n}$ ($k$ is fixed). If
there are integers $a\ge k$ and
$b\ge k$ such that the $a$- and $b$-variances of the row- and
column-representatives are  ${\cal O}(\frac {m+n}{mn})$, then 
there is a blown up matrix $\BB_{m\times n}$ such
that $\A_{m \times  n} =\BB_{m \times n} + \EEE_{m \times n}$, with  
$\| \EEE_{m \times n} \| 
={\cal O} (\sqrt{m+n} )$.
\end{thm}
\begin{pf}
 The proof  gives an explicit construction for $\BB_{m\times  n}$.
In the sequel the subscripts $m$ and $n$ will be dropped.
We shall speak in terms of microarrays (genes and conditions).

Let $\y_1 ,\dots ,\y_k \in \R^m$ and $\x_1 ,\dots ,\x_k \in \R^n$ denote the
left- and right-hand side unit-norm singular vectors corresponding to
$z_1 ,\dots ,z_k$, the 
singular values of $\A$  of  order larger than $\sqrt{m+n}$.
The $k$-dimensional representatives of the genes and conditions 
-- that are  row vectors of the
$m\times k$ matrix  $\Y = (\y_1 ,\dots ,\y_k )$ and those of the
$n\times k$ matrix  $\X = (\x_1 ,\dots ,\x_k )$, respectively --
  by the condition of the theorem
 form $a$ and $b$ clusters in $\R^k$, respectively  with sum of inner
variances  ${\cal O} ( \frac{m+n}{mn})$.
Reorder the rows and columns of $\A$ according to the clusters.
Denote by $\y^{1} ,\dots ,\y^{m} \in \R^k$ and
$\x^{1} ,\dots ,\x^{n} \in \R^k$
the Euclidean representatives
of the genes and conditions (the rows of the reordered $\Y$ and $\X$), and let
${\bar \y}^{1} , \dots , {\bar \y}^{a} \in \R^k$ and
${\bar \x}^{1} , \dots , {\bar \x}^{b} \in \R^k$
denote the cluster centers, respectively.
Now let us choose the following new representation of the genes and conditions.
The genes' representatives be row
vectors of the $m\times k$ matrix $\widetilde \Y$ such that the 
first $m_1$ rows
of $\widetilde \Y$ be equal to ${\bar \y}^{1}$, the next $m_2$ rows
to ${\bar \y}^{2}$, and so on, the  last $m_a$
rows of $\widetilde \Y$ be equal to ${\bar \y}^{a}$; similarly,
the conditions' representatives be row
vectors of the $n\times k$ matrix $\widetilde \X$ such that the first 
$n_1$ rows
of $\widetilde \X$ be equal to ${\bar \x}^{1}$, and so on, the last $n_b$
rows of $\widetilde \X$ be equal to ${\bar \x}^{b}$.

By the considerations of Theorem~\ref{thm-SA} and the assumption for the clusters,
\begin{equation}\label{dist-F-S}
 \sum_{i=1}^k {\dist}^2 (\y_i ,F ) =S_a^2 (\Y ) = {\cal O} (\frac{m+n}{mn})
\end{equation}
and
\begin{equation}\label{dist-G-S}
 \sum_{i=1}^k {\dist}^2 (\x_i ,G ) = S_b^2 (\X )= {\cal O} (\frac{m+n}{mn})
\end{equation}
hold respectively, where the $k$-dimensional subspace $F\subset \R^m$
is spanned by the column vectors of $\widetilde \Y$,
while  the $k$-dimensional subspace $G \subset \R^n$
is spanned by the column vectors of $\widetilde \X$.
We follow the construction given in~\cite{bolla1} (see Proposition 2)
 of a set  $\vv_1 ,\dots ,\vv_k $ of orthonormal vectors
within $F$ and another set $\uu_1 ,\dots ,\uu_k$ of orthonormal vectors
within $G$ such that
\begin{equation}\label{y-v}
\sum_{i=1}^k \| \y_i -\vv_i \|^2 =
\min_{\vv'_1,\dots,\vv'_k}
 \sum_{i=1}^k \| \y_i -\vv'_i \|^2
 \le 2 \sum_{i=1}^k {\dist}^2 (\y_i ,F )
\end{equation}
and
\begin{equation}\label{x-u}
 \sum_{i=1}^k \| \x_i -\uu_i \|^2 =
\min_{\uu'_1,\dots,\uu'_k}
 \sum_{i=1}^k \| \x_i -\uu'_i \|^2
 \le 2 \sum_{i=1}^k {\dist}^2 (\x_i ,G )
\end{equation}
hold,  where the minimum is taken over orthonormal sets of vectors
$\vv'_1,\dots,\vv'_k \in F$ and $\uu'_1,\dots,\uu'_k \in G$, respectively.
The construction of the vectors $\vv_1 ,\dots ,\vv_k$ is as follows 
($\uu_1 ,\dots ,\uu_k$ can be constructed
in the same way). Let $\vv'_1 ,\dots ,\vv'_k \in F$ an arbitrary
orthonormal system (obtained e.g., by the Schmidt orthogonalization method).
Let $\V' = (\vv'_1 ,\dots ,\vv'_k )$ be $m\times k$ matrix and
$$
 \Y^T \V' = \QQ \SS \Z^T
$$
be SVD, where the matrix $\SS$ contains the singular values of the $k\times k$ 
matrix $\Y^T \V'$ in its main diagonal and zeros otherwise, 
while $\QQ$ and $\Z$ are $k\times k$ orthogonal matrices
(containing the corresponding unit norm singular vector pairs in 
their columns).
The orthogonal matrix $\RRR =\Z \QQ^T$ will give the convenient orthogonal
rotation of the vectors $\vv'_1 ,\dots ,\vv'_k$.
That is, the column vectors of the matrix
$\V =\V' \RRR$ form also an orthonormal set that is the desired set
$\vv_1 ,\dots ,\vv_k$.

Define the  
 error terms $\r_i$ and $\q_i$, respectively:
$$
\r_i = \y_i - \vv_i
\quad \text{and} \quad
  \q_i = \x_i - \uu_i
  \qquad (i=1,\dots ,k).
$$
In view of (\ref{dist-F-S}) -- (\ref{x-u}), 
\begin{equation}\label{sum-r}
 \sum_{i=1}^k \| \r_i \|^2 ={\cal O} (\frac{m+n}{mn})  \quad \text{and} \quad
   \sum_{i=1}^k \| \q_i \|^2 ={\cal O} (\frac{m+n}{mn}).  
\end{equation}

Consider the following decomposition:
$$
\A = \sum_{i=1}^k z_i \y_i \x_i^T 
+\sum_{i=k+1}^{\min\{ m,n\} } z_i \y_i \x_i^T.
$$
The spectral norm of the second term is at most of order $\sqrt{m+n}$.
Now consider the first term,
\begin{equation}\label{zyx}
\begin{split}
 \sum_{i=1}^k z_i \y_i \x_i^T &=
 \sum_{i=1}^k z_i ( \vv_i +\r_i )
                        ( \uu_i^T +\q_i^T ) = \\
 &=\sum_{i=1}^k z_i  \vv_i \uu_i^T +
  \sum_{i=1}^k z_i  \vv_i \q_i^T +
   \sum_{i=1}^k z_i \r_i \uu_i^T +
 \sum_{i=1}^k z_i  \r_i \q_i^T.
\end{split}
\end{equation}
Since $\vv_1 ,\dots ,\vv_k$ and 
$\uu_1 ,\dots ,\uu_k$ are unit  vectors, 
the last three  terms in (\ref{zyx}) can be estimated by means of the relations
\begin{eqnarray*}
 \| \vv_i \uu_i^T \| &=& \sqrt{\| \vv_i \uu_i^T \uu_i \vv_i^T \| } = 1
   \qquad (i=1,\dots ,k) ,  \\
\| \vv_i \q_i^T \| &=& \sqrt{\| \q_i \vv_i^T \vv_i \q_i^T \| } =\| \q_i \|
   \qquad (i=1,\dots ,k) , \\
 \| \r_i \uu_i^T \| &=& \sqrt{\| \r_i \uu_i^T \uu_i \r_i^T \| } = \| \r_i \| 
   \qquad (i=1,\dots ,k) , \\
 \| \r_i \q_i^T \| &=&\sqrt{\| \r_i \q_i^T \q_i \r_i^T \| } =  
 \| \q_i \| \cdot  \| \r_i \|   \qquad (i=1,\dots ,k) .
\end{eqnarray*}
Taking into account that $z_i$ cannot exceed $\Theta (\sqrt{mn})$
and $k$ is fixed, 
due to (\ref{sum-r}) we get that the spectral norms of the  
last three terms in (\ref{zyx}) 
-- for their finitely many subterms  the triangle inequality is 
applicable -- are at most of order  $\sqrt{m+n}$.
Let $\BB$ be the first term, i.e., 
$$
\BB= \sum_{i=1}^k z_i \vv_i \uu_i^T   ,
$$
then $\| \A - \BB \| = {\cal O}(\sqrt{m+n})$.

By definition, the  vectors $\vv_1 ,\dots ,\vv_k$  
and the vectors $\uu_1 ,\dots ,\uu_k$ are
in the subspaces $F$ and $G$, respectively. Both spaces consist of piecewise
constant vectors,
thus the matrix $\BB$ is a blown up matrix containing $a \times b$ blocks.
The 'noise' matrix is  
$$
\EEE= \sum_{i=1}^k z_i  \vv_i \q_i^T +
   \sum_{i=1}^k z_i \r_i \uu_i^T +\sum_{i=1}^k z_i  \r_i \q_j^T+
\sum_{i=k+1}^{\min\{ m,n\} } z_i \y_i \x_i^T
$$
that finishes the proof.
\qed
\end{pf}

Then, provided the conditions of Theorem~\ref{constr} hold, by the
construction given in the proof above, an algorithm can be written 
that uses several SVD's and produces the blown up matrix $\BB$.
This $\BB$ can be regarded as the best blown up  approximation 
of the  microarray $\A$.
At the same time clusters of the genes and conditions are also
obtained. More precisely, first we conclude the clusters from the
SVD of $\A$, rearrange
the rows and columns of $\A$ accordingly, and after we use the
above construction. If we decide to perform correspondence analysis on
$\A$ then by  (\ref{A-corr}) and (\ref{DBD}), $\BB_{corr}$ will
give a good approximation
to $\A_{corr}$ and similarly, the correspondence vectors obtained by the
SVD of $\BB_{corr}$ will give representatives of the genes and
conditions.

To obtain SVD of large matrices, randomized algorithms are at our disposal,
e.g.,~\cite{AM}. There is nothing to loose when applying these
algorithms because they give
the required results only if our matrix had a primary linear structure.

\begin{ack}
The authors are indebted to G\'abor Tusn\'ady for suggesting the
microarray problem and supplying computer simulations, and also to
P\'eter Major for inspiring discussions on Chernoff-type inequalities.
\end{ack}

\end{document}